\documentclass[12pt]{article}

\usepackage{theorem}
\usepackage{amssymb}
\usepackage{amsmath}
\usepackage{amsfonts}
\usepackage{url}

\newtheorem{theorem}{Theorem}[section]
\newtheorem{lemma}[theorem]{Lemma}

\newtheorem{corollary}[theorem]{Corollary}

\def \RR {\mathbb R}
\def \EE {\mathbb E}

\def \cF {\mathcal F}
\def \cG {\mathcal G}
\def \cK {\mathcal K}
\def \cT {\mathcal T}
\def \cH {\mathcal H}

\DeclareMathOperator{\vol}{Vol}
\DeclareMathOperator{\conv}{conv}
\DeclareMathOperator{\Prob}{\mathbb P}
\DeclareMathOperator{\spn}{sp}
\DeclareMathOperator{\asp}{asp}

\def\myffrac#1#2 in #3{\raise 2.6pt\hbox{$#3 #1$}\mkern-1.5mu\raise 0.8pt\hbox{$
#3/$}\mkern-1.1mu\lower 1.5pt\hbox{$#3 #2$}}
\newcommand{\ffrac}[2]{\mathchoice%
    {\myffrac{#1}{#2} in \scriptstyle}
    {\myffrac{#1}{#2} in \scriptstyle}
    {\myffrac{#1}{#2} in \scriptscriptstyle}
    {\myffrac{#1}{#2} in \scriptscriptstyle}
}

\begin{document}

\title{On the hyperplane
conjecture for random convex sets}
\author{Bo'az Klartag\thanks{Supported by the Clay Mathematics Institute and by  NSF grant
\#DMS-0456590}  \ and \ Gady Kozma}
\date{}
 \maketitle

\begin{abstract}
Let $N \geq n+1$, and denote by $\cK$ the convex hull of $N$
independent standard gaussian random vectors in $\RR^n$. We prove
that with high probability, the isotropic constant of $\cK$ is
bounded by a universal constant. Thus we verify the hyperplane
conjecture for the class of gaussian random polytopes.
\end{abstract}

\section{Introduction}

The hyperplane conjecture suggests a positive answer to the
following question: Is there a universal constant $c
> 0$, such that for any dimension $n$ and for any convex set $\cK
\subset \RR^n$ of volume one, there exists at least one hyperplane
$\cH \subset \RR^n$ with  $\vol_{n-1}(\cK \cap \cH) > c$? Here, of
course, $\vol_{n-1}$ denotes $(n-1)$-dimensional volume.

\medskip This seemingly innocuous question, considered two decades ago by
Bourgain \cite{bou1, bou2}, has not been answered yet. We refer
the reader to, e.g., \cite{ba1}, \cite{MP} or \cite{quarter} for
partial results, history and for additional literature regarding
this hyperplane conjecture. In particular, there are large classes
of convex bodies for which an affirmative answer to the above
question is known. These include unconditional convex bodies
\cite{bou1, MP}, zonoids, duals to zonoids \cite{ba1}, bodies with
a bounded outer volume ratio \cite{MP}, unit balls of Schatten
norms \cite{KMP} and others (e.g., \cite{junge2, emanuel}).

\medskip A potential counter-example to the hyperplane
conjecture could have stemmed from random convex bodies, that
typically belong to none of these classes. Recall that, starting
with Gluskin's work \cite{gluskin}, random polytopes are a major
source of counter examples in high-dimensional  convex geometry
(in addition to the distance problem \cite{gluskin}, one has, e.g., the basis
problem \cite{szarek} or the gaussian perimeter problem \cite{nazarov}).
The goal of this short note is to show that gaussian random
polytopes, and related models of random convex sets, do not
constitute a counter-example to the hyperplane conjecture.

\medskip Suppose $\cK
\subset \RR^n$ is a convex body. The isotropic constant of $\cK$,
denoted here by $L_{\cK}$, is defined by
\begin{equation}
 n L_{\cK}^2 = \inf_{T: \RR^n \rightarrow \RR^n}
\frac{1}{\vol_n(\cK)^{1 + \frac{2}{n}}} \int_{\cK} |T x|^2 dx,
\label{def_L}
\end{equation}
 where the infimum runs over all volume-preserving
affine maps $T: \RR^n \rightarrow \RR^n$, and $| \cdot |$ stands
for the standard Euclidean norm in $\RR^n$. Directly from the
definition, the isotropic
constant is invariant under affine transformations. It is
well-known (see \cite{H} or \cite{MP}) that when $\vol_n(\cK) = 1$,
\begin{equation}
 \frac{c_1}{L_{\cK}} \leq  \inf_{T: \RR^n \rightarrow \RR^n}
\sup_{\cH \subset \RR^n} \vol_{n-1} (T\cK \cap \cH) \leq
\frac{c_2}{L_{\cK}}, \label{eq_705}
\end{equation}
 where the infimum runs
over all volume-preserving affine maps $T: \RR^n \rightarrow
\RR^n$, the supremum runs over all hyperplanes $\cH \subset
\RR^n$, and $c_1, c_2 > 0$ are some universal constants.
Throughout the text, the symbols $c, C, c^{\prime}, C^{\prime},
c_1,c_2$ etc. denote various positive universal constants, whose
value may change from one line to the next.

\medskip Thus, according to (\ref{eq_705}), the hyperplane conjecture is
equivalent to the existence of a universal upper bound for the
isotropic constant of an arbitrary convex body in an arbitrary
dimension. It is well-known that $L_{\cK} > c$ for any
finite-dimensional convex body $\cK$ (see, e.g., \cite{MP}). The
best known general upper bound is $L_{\cK} < C n^{1/4}$ for a
convex body $\cK \subset \RR^n$ (see \cite{quarter} and also
\cite{bou_dist_pols}, \cite{bou_psi_2} and \cite{dar}).

\medskip There are two natural models for random convex bodies (a body is a
compact with non-empty interior). In the centrally-symmetric
context, the first model is the symmetric convex hull of $N$
random independent points, while the second (its dual), is
the intersection of $N$ random strips. For the second model,
however, it is quite easy to demonstrate the hyperplane
conjecture. Indeed, if $N \geq 2n$, then a simple calculation
shows that it has a bounded outer volume ratio, and hence has a
bounded isotropic constant (see  \cite{MP}). If $n \leq N < 2n$,
then the hyperplane conjecture holds deterministically whenever
the resulting set is a body, according to \cite{junge}. Thus we
will focus on the first model.

\medskip We say that a random vector $X = (X_1,...,X_n) \in \RR^n$
is a standard gaussian vector if its coordinates $X_1,...,X_n$ are
independent, standard normal variables.

\begin{theorem} Let $n \geq 1, N \geq n$ and let $G_0,...,G_N$ be
independent standard gaussian vectors in $\RR^n$. Denote
$$ \cK = \conv \{ G_0,...,G_N \} \ \ \ \text{and} \ \ \ \cT = \conv \{ \pm G_1,..., \pm G_N \} $$
where $\conv$ denotes convex hull. Then, with probability greater
than $1 - C e^{-c n}$,
$$ L_{\cK} < C \ \ \ \text{and} \ \ \ L_{\cT} < C. $$
Here, $C,c > 0$ are universal constants. \label{main_thm}
\end{theorem}

Let us sketch the proof in the case that $N>n^2$ (for smaller $N$
the argument is only slightly more opaque). It is easy to see that
in this case the radius of the inscribed ball is, with high
probability, $\geq c \sqrt{\log N}$ --- just calculate the
probability that in a given direction all points are inside a
strip of width $\epsilon \sqrt{\log N}$ and do a union-bound over
a dense net in $S^{n-1}$. On the other hand, with probability $1$
all faces of $\cK$ are $(n-1)$-dimensional simplices. Further,
with high probability all centers of gravity of all simplices are
with distance $\leq C\sqrt{\log N}$ from $0$ --- the center of
gravity of each simplex is a gaussian vector whose coordinates
have variance $\ffrac{1}{n}$, and again all that is needed is to
do a union-bound over all $n$-tuples of vertices. The
concentration of the volume of a simplex around its center of
gravity shows that almost all the mass is within distance of
$\approx \sqrt{\log N}$ which implies the required estimate for
$L_\cK$.

\medskip On first glance it seems that some miracle is at work --- why should
we get the same $\sqrt{\log N}$ in the lower and upper bound?
However, this is some manifestation of the phenomenon that the
maximum of many independent variables is strongly concentrated. Of
course, the different faces of our body are not independent, but
it turns out that they are sufficiently independent to display a
similar concentration phenomenon. Thus our proof is robust and
would admit direct generalizations to other types of
distributions, in place of the standard gaussian distribution. We
will prove it for some other distributions that include the
uniform distribution on the cube and on its corners $\{\pm 1\}^n$,
see Theorem \ref{thm_1149} below. The technique should also work
for points uniform on $S^{n-1}$.

\medskip
\emph{Acknowledgement.} We would like to thank Jean Bourgain for
motivating us to work on this problem. BK would
also like to thank Alain Pajor for discussions on the subject.

\section{Simplices}
\label{section2}

In this section, we assume that $N \geq n \geq 1$ are integers,
and that $G_0,...,G_N$ are independent random vectors in $\RR^n$ which need
not be identically distributed.
We write $G_i = (G_{i,1},...,G_{i, n}) \in \RR^n$, and we make the
following assumptions regarding $G_0,...,G_N$:

\begin{enumerate}
\item[($\ast$a)] The random variables $G_{i,j} \ (i=0,...,N, \, j=1,...,n)$ are independent.
\item[($\ast$b)] For any $i=0,...,N, \, j=1,...,n$,
$$ \EE G_{i,j} =0, \ \EE G_{i,j}^2 = 1 \ \ \text{and} \ \ \EE \exp \left( G_{i,j}^2 /
10 \right) \leq 10. $$
\item[($\ast$c)] The  $G_i \ (i=0,...,N)$
  are absolutely continuous.
\end{enumerate}

The constant $10$ plays no special r\^ole. Note that ($\ast$a),
($\ast$b) and ($\ast$c) hold when $G_0,...,G_N$ are independent
standard gaussian vectors. Our main technical tool is the
following Bernstein's inequality for variables with exponential
tail (``$\psi_1$''), see, e.g.\ \cite[Section 2.2.2]{vw}.

\begin{theorem}
Let $L > 0$, let $m \geq 1$ be an integer, and let $X_1,...,X_m$
be independent random variables with zero mean. Assume that $$ \EE
\exp (|X_i| / L) \leq 20 \ \ \ \text{for} \ \ 1 \leq i \leq m. $$
Then, for any $t
> 0$,
$$
\Prob \left \{ \left|\frac{1}{m} \sum_{i=1}^m X_i \right|  > t
\right \} \leq 2 \exp \left( -c m \min \left \{ \frac{t}{L},
\frac{t^2}{L^2} \right \} \right),
$$
where $c > 0$ is a universal constant. \label{bernstein}
\end{theorem}

 The following lemma is a consequence of Theorem \ref{bernstein}.

\begin{lemma} Fix $0 \leq k_1 < k_2 < \dotsb < k_n \leq N$. Then,
\begin{enumerate}
\item[(i)] $\displaystyle \Prob \left \{ \left| \frac{1}{n} \sum_{i=1}^n
G_{k_i} \right| > C \sqrt {\log \frac{2 N}{n}}  \right \} < \left(
\frac{n}{10 N} \right)^n,$
\item[(ii)] $\displaystyle \Prob \left \{ \left|\frac{1}{N+1} \sum_{i=0}^N
G_{i} \right| > C \sqrt {\log \frac{2 N}{n}}  \right \} < \left(
\frac{n}{10 N} \right)^n,$
\item[(iii)] $\displaystyle \Prob \left \{
\frac{1}{n^2} \sum_{i=1}^n \left[ \left( \sum_{j=1}^n G_{k_i,j}
\right)^2 + \sum_{j=1}^n G_{k_i,j}^2 \right]
> C \log \frac{2 N}{n}  \right \} < \left(
\frac{n}{10 N} \right)^n.$
\end{enumerate}
Here, $C > 0$ is a universal constant. \label{lem_845}
\end{lemma}

\noindent\emph{Proof.} Fix $1 \leq i, j \leq n$. We have, for any
$t$,
$$
\EE\exp(tG_{k_i,j})\leq
\exp(\tfrac{5}{2}t^2)\EE\exp(G_{k_i,j}^2/10) \stackrel{(\ast
\text{b})}{\leq} 10\exp(\tfrac{5}{2}t^2).
$$
By independence, for any $j=1,..,n$ and $t$,
$$
\EE \exp \left( t \sum_{i=1}^n G_{k_i,j} \right) \leq 10 \exp
\left( \tfrac{5}{2} n t^2 \right),
$$
so
\begin{equation}
\EE\exp  \left( \left|t\sum_{i=1}^n G_{k_i,j}\right| \right) \leq
\EE \exp \left( t \sum_{i=1}^n G_{k_i,j} \right) +
    \exp \left( -t \sum_{i=1}^n G_{k_i,j} \right) \leq
20\exp(\tfrac{5}{2}nt^2). \label{eq_946}
\end{equation}
 Denote $Y_j = \frac{1}{\sqrt{n}}
\sum_{i=1}^n G_{k_i,j} \ \ (j=1,...,n)$. Then the $Y_j$ are
independent random variables of mean zero and variance one.
Moreover, by (\ref{eq_946}), for $j=1,...,n$,
\begin{align}
\nonumber \EE \exp \left( \frac{|Y_j^2 - 1|}{100} \right) &\leq 2
\EE \exp \left( \frac{Y_j^2}{100} \right)   2 \int_0^{\infty} \Prob \left\{\frac{Y_j^2}{100} > \log t \right\} \,dt\leq \\
&\leq 2 + 2 \int_1^{\infty} t^{-10} \EE e^{\sqrt{\log t} |Y_j|}
\,dt
  \stackrel{(\ref{eq_946})}{\leq}
  2 + 40\int_1^{\infty} t^{-15/2} \,dt < 9.
\label{eq_139}
\end{align}
Hence, we may apply Theorem \ref{bernstein} for the independent
random variables $Y_j^2 - 1$, with $L = 100$. We conclude that
\begin{multline*}
\Prob \left \{ \left| \frac{1}{n} \sum_{i=1}^n G_{k_i} \right|
> C \sqrt{\log \frac{2 N}{n}} \right \}  \Prob \left \{  \frac{1}{n} \sum_{j=1}^n Y_j^2
> C^2 \log \frac{2 N}{n} \right \} \leq \\ \leq
\Prob \left \{  \frac{1}{n} \sum_{j=1}^n [Y_j^2 - 1]
> (C/2)^2 \log \frac{2 N}{n} \right \}
\leq 2 \exp \left( -c n \frac{C^2}{400} \cdot \log \frac{2 N}{n}
\right) < \left( \frac{n}{10 N} \right)^n,
\end{multline*}
for an appropriate choice of a large universal constant $C > 50$.
This completes the proof of (i). The argument that leads to (ii)
is similar. We define $\tilde{Y}_j = \frac{1}{\sqrt{N+1}}
\sum_{i=0}^N G_{i,j} \ \ (j=1,...,n)$. Arguing exactly as above,
we find  that for $j=1,...,n$,
$$ \EE \exp \left( \frac{|\tilde{Y}_j^2 - 1|}{100} \right) < 9.
$$
We may invoke Theorem \ref{bernstein} for the independent random
variables $\tilde{Y}_j^2 - 1$, with $L = 100$. This yields
\begin{multline*}
\Prob \left \{ \left| \frac{1}{N+1} \sum_{i=0}^N G_{i} \right|
> C \sqrt{\log \frac{2 N}{n}} \right \}
= \Prob \left \{  \frac{1}{N+1} \sum_{j=1}^n \tilde{Y}_j^2
> C^2 \log \frac{2 N}{n} \right \} \leq \\ \leq
\Prob \left \{  \frac{1}{n} \sum_{j=1}^n [\tilde{Y}_j^2 - 1]
> (C/2)^2 \log \frac{2 N}{n} \right \}
\leq 2 \exp \left( -c n \frac{C^2}{400} \cdot \log \frac{2 N}{n}
\right) < \left( \frac{n}{10 N} \right)^n,
\end{multline*}
provided that $C$ is a sufficiently large universal constant. This
proves (ii). It remains to prove (iii). The random variables
$G_{i,j}^2 - 1 \ (i=0,...,N, \ j=1,...,n)$ are independent, have
mean zero and they satisfy that
$$
 \EE \exp \left[ |G_{i,j}^2 -1| / 10 \right] \leq 2 \EE \exp
(G_{i,j}^2 / 10) \leq 20, $$ according to ($\ast$b). Hence, we may
apply Theorem \ref{bernstein} for the independent random variables
$G_{k_i,j}^2 - 1$, with $L = 10$. We conclude that for any $1 \leq
j \leq n$,
\begin{equation}
 \Prob \left \{  \frac{1}{n} \sum_{i=1}^n [G_{k_i, j}^2 - 1] > C
\log \frac{2 N}{n} \right \} \leq 2 \exp \left( -c n \frac{C}{10}
\cdot \log \frac{2 N}{n} \right) < \left( \frac{n}{20 N}
\right)^n, \label{eq_830}
\end{equation}
for a large universal constant $C$. We sum (\ref{eq_830}) over
$j=1,...,n$ and conclude that
\begin{eqnarray}
\label{eq_931} \lefteqn{ \Prob \left \{ \frac{1}{n^2}
\sum_{i,j=1}^n G_{k_i,j}^2
> (C +2) \log \frac{2 N}{n} \right \} } \\ & \leq & \sum_{j=1}^n \Prob \left \{
\frac{1}{n} \sum_{i=1}^n [G_{k_i, j}^2 - 1] > C \log \frac{2 N}{n}
\right \} \stackrel{(\ref{eq_830})}{<} n \left( \frac{n}{20 N}
\right)^n < \frac{1}{2} \cdot \left( \frac{n}{10 N} \right)^n.
\nonumber
\end{eqnarray}
Denote $Z_i = \frac{1}{\sqrt{n}} \sum_{j=1}^n G_{k_i,j} \ \
(i=1,...,n)$. Then the $Z_i$ are independent and $\EE Z_i^2 = 1$
for $i=1,...,n$. Repeating the argument from (\ref{eq_946}) and
(\ref{eq_139}) (the only difference is that here we sum over  $j$ and before we
summed over $i$), we obtain that $ \EE \exp \left( |Z_i^2 - 1| / 100
\right) < 9$. Thus, we may use Theorem \ref{bernstein} for the
random variables $Z_i^2 - 1$ with $L = 100$, and deduce that
\begin{equation}
\Prob \left \{ \frac{1}{n} \sum_{i=1}^n [Z_i^2 - 1]
> C \log \frac{2 N}{n}
 \right \} \leq 2 \exp \left( -c n \frac{C}{100} \log \frac{2 N}{n}
 \right) < \frac{1}{2} \cdot \left( \frac{n}{10 N} \right)^n.
 \label{eq_1002}
\end{equation}
The desired conclusion (iii) follows at once from (\ref{eq_931})
and (\ref{eq_1002}). \hfill $\square$

\medskip
The next lemma is a simple, concrete calculation for the regular
$(n-1)$-simplex. We write $e_1,...,e_n$ for the standard basis in $\RR^n$, and
denote $$ \triangle^{n-1} = \conv \{e_1,...,e_n \} \subset
\RR^n, $$ the $(n-1)$-dimensional regular simplex.

\begin{lemma} Let $X = (X_1,...,X_n)$
be a random vector that is distributed uniformly in
$\triangle^{n-1}$. Then,
$$ \EE X_i X_j = \frac{1 + \delta_{i,j}}{n(n+1)} $$
where $\delta_{i,j}$ is Kronecker's delta. \label{lem_906}
\end{lemma}

\noindent\emph{Proof.} Examine $X$ without its last coordinate,
$(X_1,...,X_{n-1})$. This is distributed
uniformly in the simplex $$ \left \{ x \in \RR^{n-1} \ ; \
\sum_{i=1}^{n-1} x_i \leq 1, \ \forall i, x_i \geq 0 \right \}. $$
Consequently, the density of the random variable $Y = X_1 + \dotsb
+ X_{n-1}$ is proportional to $t \mapsto t^{n-2}$ in the interval
$(0,1)$, and is zero elsewhere. Hence,
\begin{equation}
\EE (X_1 + \dotsb + X_{n-1})^2 = \EE Y^2 = \int_0^1 t^2 \cdot (n-1)
t^{n-2} dt = \frac{n-1}{n+1}. \label{eq_1042}
\end{equation}
Note that $\sum_{i=1}^n X_i \equiv 1$. Therefore
\begin{equation}
1 = \EE \left( \sum_{i=1}^n X_i \right)^2 =  \sum_{i=1}^n \EE
X_i^2 + \sum_{i \neq j} \EE X_i X_j. \label{eq_1043}
\end{equation}
From (\ref{eq_1042}) and (\ref{eq_1043}) we get that, when $n \geq
2$,
$$ (n-1) \EE X_1^2 + (n-1) (n-2) \EE
X_1 X_2 = \frac{n-1}{n+1},  \ \ \ n \EE X_1^2 + n (n-1) \EE X_1
X_2 = 1, $$ and  the lemma follows. \hfill $\square$

\begin{corollary} Fix $0 \leq k_1 < k_2 < ... < k_n \leq N$,
and set $\cF = \conv \{ G_{k_1},...,G_{k_n} \}$. Denote $Z \frac{1}{N+1} \sum_{i=0}^N G_i$. Then with probability greater
than $1 - 4 \left( \frac{n}{10 N} \right)^n$, the set $\cF$ is an
$(n-1)$-dimensional simplex that satisfies
\begin{enumerate}
\item[(i)] $\displaystyle \frac{1}{\vol_{n-1}(\cF)} \int_{\cF} |x|^2 dx <  C \log \frac{2
N}{n}, $
\item[(ii)] $\displaystyle  \frac{1}{\vol_{n-1}(\cF)} \int_{\cF} |x - Z|^2 dx <  C \log \frac{2
N}{n} $.
\end{enumerate}
Here, $C > 0$ is a universal constant. \label{cor_1058}
\end{corollary}

\noindent\emph{Proof.}  The random vectors $G_{k_i}$ are
independent and absolutely continuous according to ($\ast$c).
Hence, with probability one, the linear span of
$G_{k_1},...,G_{k_n}$ equals $\RR^n$, and the set $\cF$ is an
$(n-1)$-dimensional simplex in $\RR^n$ whose vertices are the
points $G_{k_1},...,G_{k_n}$. Denote $T (G_{k_i,j})_{i,j=1,...,n}$, an $n \times n$ matrix. Then, $ \cF T \left( \triangle^{n-1} \right) $ and hence,
\begin{equation}
 \frac{1}{\vol_{n-1}(\cF)} \int_{\cF} |x|^2 dx \frac{1}{\vol_{n-1}(\triangle^{n-1})} \int_{\triangle^{n-1}} |Tx|^2\,
dx. \label{eq_1039} \end{equation} According to Lemma
\ref{lem_906},
\begin{eqnarray}
\label{eq_1040} \lefteqn{ \frac{1}{\vol_{n-1}(\triangle^{n-1})}
\int_{\triangle^{n-1}} |Tx|^2 \,dx = \sum_{i=1}^n \sum_{j_1, j_2
=1}^n G_{k_i, j_1} G_{k_i, j_2}
 \frac{1}{\vol_{n-1}(\triangle^{n-1})} \int_{\triangle^{n-1}}
 x_{j_1}
 x_{j_2} \,dx } \\ & = & \frac{1}{n(n+1)} \sum_{i=1}^n \sum_{j_1, j_2
=1}^n G_{k_i, j_1} G_{k_i, j_2} (1 + \delta_{j_1,j_2}) \frac{1}{n(n+1)} \sum_{i=1}^n \left[ \left( \sum_{j=1}^n G_{k_i,j}
\right)^2 + \sum_{j=1}^n G_{k_i,j}^2 \right]. \nonumber
 \end{eqnarray}
We conclude from (\ref{eq_1039}), (\ref{eq_1040}) and Lemma
\ref{lem_845}(iii) that
\begin{equation}
\Prob \left \{ \frac{1}{\vol_{n-1}(\cF)} \int_{\cF} |x|^2 dx > C
\log \frac{2 N}{n}  \right \} < \left( \frac{n}{10 N} \right)^n.
\label{eq_1042_}
\end{equation}

As for the second part of the corollary, according to Lemma \ref{lem_845}(i)
and Lemma \ref{lem_845}(ii) we know that
\begin{equation}
\Prob \left \{  \left| \frac{1}{n} \sum_{i=1}^n G_{k_i} \right| +
\left| \frac{1}{N+1} \sum_{j=0}^N G_j \right| < 2 C \sqrt{\log
\frac{2 N}{n}} \right \} \geq 1 - 2 \left( \frac{n}{10 N}
\right)^n. \label{eq_848}
\end{equation}
Additionally,
$$ \frac{1}{\vol_{n-1}(\cF)} \int_{\cF} |x - Z|^2 \,dx = |Z|^2 - 2 \left
\langle \frac{1}{\vol_{n-1}(\cF)} \int_{\cF} x\, dx, Z \right \rangle
+ \frac{1}{\vol_{n-1}(\cF)} \int_{\cF} |x|^2 \,dx  $$
\begin{equation}  = \left| \frac{1}{N+1}
\sum_{j=0}^N G_j \right|^2 - 2 \left \langle \frac{1}{n} \sum_{i=1}^n
G_{k_i},
\frac{1}{N+1}
\sum_{j=0}^N G_j \right \rangle + \frac{1}{\vol_{n-1}(\cF)} \int_{\cF}
|x|^2 \,dx. \label{eq_1048}
\end{equation}
 By combining (\ref{eq_1048}) with
(\ref{eq_1042_}) and (\ref{eq_848}), we obtain
\begin{equation}
\Prob \left \{ \frac{1}{\vol_{n-1}(\cF)} \int_{\cF} |x - Z|^2 \,dx > C
\log \frac{2N}{n} \right \} < 3 \left( \frac{n}{10 N} \right)^n.
\label{eq_1047}
\end{equation}
From (\ref{eq_1042_}) and (\ref{eq_1047}) the corollary follows.
\hfill $\square$

\medskip For a point $x \in \RR^n$ and a set $A \subset \RR^n$, we write
$d(x, A) = \inf_{y \in A} |x-y|$.

\begin{lemma} Set $$ \cK = \conv \{ G_{0},...,G_{N} \}, \  \ \
\cT = \conv \{ \pm G_1,..., \pm G_N \}, $$ and denote $Z \frac{1}{N+1} \sum_{i=0}^N G_i$. Then with probability greater
than $1 - C e^{-c n}$,
\begin{enumerate}
\item[(i)] $ \displaystyle
\frac{1}{\vol_n(\cT)} \int_{\cT} |x|^2 dx < C \log \frac{2 N}{n},$
\item[(ii)] $ \displaystyle
\frac{1}{\vol_n(\cK)} \int_{\cK} |x - Z|^2 dx < C \log \frac{2
N}{n}.$
\end{enumerate}
Here, $c, C > 0$ are universal constants. \label{cor_1059}
\end{lemma}

\noindent\emph{Proof.} The random vectors $G_{i,j}$ are absolutely
continuous, by assumption ($\ast$c). Hence, with probability one,
the points $ G_0,...,G_N$ are in general position in $\RR^n$; that
is, with probability one, no $n+1$ distinct points from $\{
G_0,...,G_N \}$ lie in the same affine hyperplane in $\RR^n$.
Consequently, all the $(n-1)$-dimensional facets of the polytopes
$\cK$ and $\cT$ are simplices, with probability one. Note that, in
the case of $\cT$ we use the fact that $G_i$ and $-G_i$ could
never belong to the same face.

\medskip Let $\cF_1,...,\cF_\ell$ be a complete list of the
$(n-1)$-dimensional
facets of $\cT$. Since a facet is determined by $n$ points
from $\{ \pm G_1,...,\pm G_N \}$, then
$$
\ell \leq {2N \choose n} \leq \left(\frac{2eN}{n}\right)^n.
$$ According to Corollary
\ref{cor_1058}(i), with probability greater than
 $1 - 4 \left( \frac{2e}{10} \right)^n$,
\begin{equation}
 \int_{\cF_i} |x|^2  dx <  C \left(\log
\frac{2 N}{n} \right) \cdot \vol_{n-1}(\cF_i), \ \ \ \text{for} \ i=1,...,\ell.
 \label{eq_1104}
\end{equation}
Each point $x \in \cT$
(except for the origin) may be uniquely represented as $x = t y$
with $0 < t \leq 1$ and $y \in \partial \cT$. We integrate with
respect to these standard polar coordinates, and obtain that
\begin{equation}
 \frac{1}{\vol_n(\cT)} \int_{\cT} |x|^2 \,dx = \frac{1}{\vol_n(\cT)}
\int_{\partial \cT} \int_0^1  |t y|^2 t^{n-1} \langle y, \nu_y
\rangle \,dt \,dy \label{eq_1103}
\end{equation}
where $\nu_y$ is the unit outward normal to $\partial \cT$ at $y$
($\nu_y$ is uniquely defined almost everywhere as $\cT$ is
convex). When $y \in \cF_i$ for some $i=1,...,\ell$, we have that
$\langle y, \nu_y \rangle = d(0, \asp \cF_i)$ where $\asp \cF_i$ is the affine
subspace spanned by $\cF_i$. Hence, from
(\ref{eq_1103}), \begin{equation} \frac{1}{\vol_n(\cT)} \int_{\cT}
|x|^2 dx = \frac{1}{\vol_n(\cT)} \sum_{i=1}^{\ell}
  \frac{d(0, \asp\cF_i)}{n+2} \int_{\cF_i} |y|^2
dy. \label{eq_1109}
\end{equation}
 Recall that $\sum_{i=1}^{\ell} d(0, \asp\cF_i) \cdot
\vol_{n-1}(\cF_i) = n \vol_n(\cT)$. We combine (\ref{eq_1109}) with
(\ref{eq_1104}), and conclude that with probability greater than
 $1 - 4 \left( \frac{2e }{10 } \right)^n$,
$$ \frac{1}{\vol_n(\cT)} \int_{\cT}
|x|^2 dx  <  \frac{1}{\vol_n(\cT)} \sum_{i=1}^{\ell} \frac{d(0,
\asp\cF_i) \cdot \vol_{n-1}(\cF_i)}{n+2} \cdot C \log \frac{2N}{n} <
C \log \frac{2N}{n}. $$ This completes the proof of (i).
The proof of (ii) is very similar; we supply some details.
Let $\cG_1,...,\cG_{k}$
denote the $(n-1)$-dimensional facets of $\cK$.
Observe that $Z \in \cK$, and that any $x \in \cK$
(except for the point $Z$) is uniquely represented as $x = Z + t (y - Z)$ with
$0 < t \leq 1$ and $y \in
\partial \cK$. As before, integration in polar coordinates yields
\begin{multline*}
\frac{1}{\vol_n(\cK)} \int_{\cK} |x - Z|^2 \,dx    \frac{1}{ \vol_n(\cK)} \int_{\partial \cK} \int_0^1 t^{n-1} |t
(y - Z) |^2 \langle y - Z, \nu_y \rangle \,dt \,dy = \\
=\frac{1}{ n+2 } \sum_{i=1}^k \frac{d(Z, \asp \cG_i)}{\vol_n(\cK)}
\int_{\cG_i} |y - Z|^2 \,dy.
\end{multline*}
Again, $\sum_{i=1}^k d(Z, \asp \cG_i) \vol_{n-1}(\cG_i) = n
\vol_n(\cK)$. Thus, in order to prove (ii), we may simply
reproduce the argument from the proof of (i), with Corollary
\ref{cor_1058}(ii) replacing the role of Corollary
\ref{cor_1058}(i). This completes the proof. \hfill $\square$

\section{Random Polytopes}
\label{section3}

We summarize the results of Section \ref{section2} in the
following corollary. Note that the convex bodies
discussed in this corollary have diameter that is larger than $c \sqrt{n}$ with high
 probability. Nevertheless, it is still possible to prove a much
better estimate regarding the second moment of the Euclidean norm.

\begin{corollary} Let $N \geq n \geq 1$ and suppose
that $G_0,...,G_N$ are independent random vectors in $\RR^n$ that
satisfy conditions ($\ast$a) and ($\ast$b) above. Set $\cK = \conv
\{ G_{0},...,G_{N} \}, \cT = \conv \{ \pm G_1,...,\pm G_N \}$,
and denote $Z = \frac{1}{N+1} \sum_{i=0}^N
G_i$. Then with probability greater than $1 - C e^{-c n}$,
$$
 \frac{1}{\vol_n(\cT)} \int_{\cT} |x|^2 dx < C \log \frac{2
N}{n} \ \ \ \text{and} \ \ \
 \frac{1}{\vol_n(\cK)} \int_{\cK} |x - Z|^2 dx < C \log \frac{2
N}{n},
$$
 where $C,c > 0$ are universal constants. \label{cor_455}
\end{corollary}

\noindent\emph{Proof.} Suppose first that $G_0,...,G_N$ satisfy also
($\ast$c); that is, assume that they are absolutely continuous random
 variables. Then the desired conclusion
 follows from Lemma \ref{cor_1059}. For the general case,
 note that the quantity
\begin{equation} \Prob \left \{
 \frac{1}{\vol_n(\cT)} \int_{\cT} |x|^2 dx < C \log \frac{2
N}{n}, \ \
\frac{1}{\vol_n(\cK)} \int_{\cK} |x - Z|^2 dx < C \log \frac{2
N}{n} \right \} \label{eq:weakcont}
\end{equation} depends continuously on the distribution of
$G_0,...,G_N$ in the weak topology at measures where
$\Prob(\vol_n(\cK)=0)=\Prob(\vol_n(\cT)=0)=0$. At other measures
(\ref{eq:weakcont}) may have a discontinuity of no more than
$\Prob(\vol_n(\cK)=0) + \Prob(\vol_n(\cT)=0)$. The corollary
follows by approximating $G_0,...,G_N$ with absolutely continuous
random vectors that satisfy ($\ast$a), ($\ast$b) and ($\ast$c) and
noting that by ($\ast$a) and ($\ast$b) we have
$\Prob(\vol_n(\cK)=0) + \Prob(\vol_n(\cT)=0) <Ce^{-cn}$, according
to \cite{rudelson}. \hfill $\square$

\medskip The next theorem is concerned with non-gaussian analogs of
Theorem \ref{main_thm}. The main new case covered by that theorem
is that of random sign vectors, i.e., independent random vectors
whose coordinates are independent, symmetric Bernoulli variables.
We remark in passing that in the Bernoulli case the probability of
$\cK$ or $\cT$ to be degenerate was known before \cite{rudelson}.
See \cite{KKS95, TV05}.

\begin{theorem} Let $n \geq 1$ and $2n \leq N \leq 2^n$. Suppose
that $G_1,...,G_N$ are independent random vectors in $\RR^n$ that
satisfy conditions ($\ast$a) and ($\ast$b) above. Set $\cT = \conv
\{ \pm G_1,...,\pm G_N \}$. Then with probability greater than $1
- Ce^{-c n}$,
$$ L_{\cT} < C $$
 where $C,c > 0$ are universal constants. \label{thm_1149}
\end{theorem}

\noindent\emph{Proof.} We may clearly assume that $n$ exceeds a certain
universal constant. It was proved in \cite[Theorem 4.8]{LPRT} that, under the
assumptions of the present theorem,
\begin{equation}
 \left( \vol_n(\cT) \right)^{1/n} > C \sqrt{ \frac{\log
(2 N/n)}{n}} \label{eq_1154} \end{equation} with probability
greater than $1- C e^{-c n}$. From Corollary \ref{cor_455} we know
that with probability larger than $1 - Ce^{-c^{\prime} n}$,
\begin{equation}
\frac{1}{\vol_n(\cT)} \int_{\cT} |x|^2 dx < C \log \frac{2 N}{n}
\label{eq_1156}
\end{equation}
The theorem follows by substituting the estimates (\ref{eq_1154})
and (\ref{eq_1156}) into the definition (\ref{def_L}). \hfill
$\square$

\medskip Generally speaking, the restrictions on $N$ in Theorem \ref{thm_1149} are typically quite
easy to work around. For example, in the case of Bernoulli
variables, if $2^n\leq N<3^n$ then (\ref{eq_1154}) and hence the
conclusion of Theorem \ref{thm_1149} hold with a different
constant, while if $N\geq 3^n$ then with very high probability
$\cT$ is a hypercube. As for $n\leq N<2n$, one can show
(\ref{eq_1154}) by noting that even a single simplex has enough
volume --- see e.g., \cite{TV06}.

\medskip Our next lemma shows the same volume estimates in the Gaussian
case. It is standard and well-known.

\begin{lemma} Let $N \geq n \geq 1$ and suppose that $G_0,...,G_N$ are
independent standard gaussian vectors in $\RR^n$. Denote $$ \cK \conv \{ G_{0},...,G_{N} \} \ \ \ \text{and} \ \ \ \cT = \conv \{
\pm G_1,...,\pm G_N \}. $$ Then, with probability greater than $1-
C e^{-c n}$,
$$ \vol_n(\cK)^{\frac{1}{n}} > c \sqrt{ \frac{\log (2N / n)}{n} }
\ \ \ \text{and} \ \ \ \vol_n(\cT)^{\frac{1}{n}} > c
\sqrt{\frac{\log (2 N / n)}{n}}, $$ where $c, C > 0$ are universal
constants. \label{known}
\end{lemma}

\noindent\emph{Proof sketch.} We start with the lower bound for
$\vol_n(\cT)$. For the range $N \geq 2n$, it is well-known (see,
e.g., \cite{Glu2,LPRT} and references therein) that with
probability greater than $1 - C e^{-c n}$,
\begin{equation}
 c  \sqrt{\log \frac{2 N}{n}} D^n \subseteq \cT
\label{eq_323}
\end{equation}
where $D^n = \{ x \in \RR^n ; |x| \leq 1 \}$ is the unit Euclidean
ball in $\RR^n$. Since $\vol^{1/n}(D^n) > c / \sqrt{n}$, the desired
lower
bound for $\vol_n(\cT)$ follows from (\ref{eq_323}) in this case.
It remains to deal
with the range $n \leq N < 2n$. It turns out that in this range a single
simplex supplies enough volume for our needs. We thus assume that $N = n$. Elementary Euclidean geometry
shows that
$$ \vol_n(\cT) = \frac{2^n}{n!} \prod_{i=1}^n d \left( G_i, \spn \{
G_1,...,G_{i-1} \} \right) $$ where $\spn$ stands for linear span
(we define $\spn (\emptyset) = \{ 0 \}$). Denote $\xi_i = d \left( G_i,
\spn \{ G_1,...,G_{i-1} \} \right)$ for $i=1,...,n$. Then $\xi_1,...,\xi_n$ are
independent random variables, with $\xi_i^2$ being distributed
chi-square with $i$ degrees of freedom. Standard estimates (one
may use, e.g., Bernstein's inequality above) show that
$$
\Prob \left \{ \vol_n(\cT) < \frac{(2c)^n}{\sqrt{n!}} \right \} \Prob \left \{  \sum_{i=1}^n \left[ -\log \frac{\xi_i}{\sqrt{i}}
\right] > n \log \frac{1}{c} \right \} < C e^{-c^{\prime} n}
$$
for a suitable choice of universal constants $c, C, c^{\prime} >
0$. This completes the proof of the desired lower bound for
$\vol_n(\cT)$. Regarding $\vol_n(\cK)$, denote $G_i^{\prime} = G_i
- G_0$ for $i=1,...,n$, and set $\cK^{\prime} = \conv \{ \pm
G_1^{\prime},...,\pm G_n^{\prime} \}$. Then,
\begin{equation}
 \vol_n (\cK) = \vol_n ( \conv \{ 0, G_1 - G_0, G_2 - G_0,..., G_N -
G_0 \} ) \geq 4^{-n} \vol_n( \cK^{\prime} ) \label{eq_1144}
\end{equation}
 by the
Rogers-Shephard inequality \cite{rs}. With probability one, the
vectors $G_1,...,G_n$ are linearly independent. Let $S: \RR^n
\rightarrow \RR^n$ be the unique linear map that satisfies $S(G_i)
= G_i - G_0$ for $i=1,...,n$. Then $\cK^{\prime} = S(\cT)$. Hence,
\begin{equation}
 \vol_n(\cK^{\prime}) = \det(S) \cdot \vol_n(\cT).
 \label{eq_111}
\end{equation}
Let $v \in \RR^n$ be such that $\langle v, G_i \rangle = 1$ for
$i=1,...,n$. The vector $v$ is independent of $G_0$. Moreover,
$|v| \geq 1 / |G_1|$, and with probability greater than $1 - C
e^{-c n}$ we have $|G_1| \leq C \sqrt{n}$. Consequently,
$$ \Prob \{ |v| < c / \sqrt{n} \} \leq C e^{-c n}. $$ Clearly, $S
x = x - \langle x, v\rangle G_0$ for all $x \in \RR^n$. Therefore
$\det(S) = 1 - \langle v, G_0 \rangle$. Conditioning on $v$, we
see that $\det(S)$ is a gaussian random variable with mean $1$ and
variance $|v|^2$. Hence,
\begin{equation}
 \Prob \{ |\det(S)| < 2^{-n} \} = \EE_v \int_{-2^{-n}}^{2^{-n}}
\frac{1}{\sqrt{2 \pi |v|^2}} \exp \left(-\frac{(t-1)^2}{2 |v|^2}
\right) dt \leq C e^{-c n} + C \sqrt{n} 2^{-n}. \label{eq_112}
\end{equation}
The desired lower bound for $\vol_n(\cK)$ follows from
(\ref{eq_1144}), (\ref{eq_111}), (\ref{eq_112}) and from the lower
bound for $\vol_n(\cT)$, that was already proven. \hfill $\square$

\medskip
\noindent\emph{Proof of Theorem \ref{main_thm}.} From Corollary
\ref{cor_455} and Lemma \ref{known}, we know that with probability
greater than $1 - C e^{-c n}$,
$$
\frac{1}{\vol_n(\cT)^{1 + \frac{2}{n}}} \int_{\cT} |x|^2 dx < C n \
\ \ \text{and} \ \ \ \frac{1}{\vol_n(\cK)^{1 + \frac{2}{n}}}
\int_{\cK} |x - Z|^2 dx < C n
$$
for some point $Z \in \RR^n$ depending on $\cK$. The theorem
follows from the definition (\ref{def_L}). \hfill $\square$




{\small  \noindent Department of Mathematics, Princeton
university, Princeton, NJ 08544, USA \\ {\it e-mail address:}
\verb"bklartag@princeton.edu"

\bigskip \noindent Department of Mathematics, Weizmann Institute of
Science, Rehovot 76100, Israel \\  {\it e-mail address:}
\verb"gady.kozma@weizmann.ac.il" }

\end{document}